\documentclass[12pt,a4paper,reqno,final]{amsart}

\usepackage{amssymb,amsmath,amsfonts,amsthm,mathrsfs}
\usepackage{xcolor}

\usepackage[margin=3cm]{geometry}

\newcommand{\N}{\mathbb N}

\newcommand{\R}{\mathbb R}
\newcommand{\Rone}{\mathbb R}

\newcommand{\la}{\left\langle}
\newcommand{\ra}{\right\rangle}

\newcommand{\mc}{\mathcal}
\newcommand{\lint}{\int\limits}

\newcommand{\lara}[1]{\langle #1 \rangle}
\newcommand{\irm}{{\rm i}}
\newcommand{\lsum}{\sum\limits}

\DeclareMathOperator{\WF}{WF}

\numberwithin{equation}{section}
\newtheorem{thm}{Theorem}[section]
\newtheorem{cor}[thm]{Corollary}

\newtheorem{remark}[thm]{Remark}
\newtheorem{definition}[thm]{Definition}

\title[Hyperbolic systems with non-diagonalisable principal part]{
	Hyperbolic systems with non-diagonalisable principal part and variable multiplicities, II: microlocal analysis}

\author[Claudia Garetto]{Claudia Garetto}
\address{
	Claudia Garetto:
	\endgraf
	Department of Mathematical Sciences
	\endgraf
	Loughborough University
	\endgraf
	Loughborough, Leicestershire, LE11 3TU
	\endgraf
	United Kingdom
	\endgraf
	{\it E-mail address} {\rm c.garetto@lboro.ac.uk}
}
\author[Christian J\"ah]{Christian J\"ah}
\address{
	Christian J\"ah:
	\endgraf
	Department of Mathematics
	\endgraf
	G\"ottingen University
	\endgraf
	G\"ottingen, D-37073
	\endgraf
	Germany
	\endgraf
	{\it E-mail address} {\rm christian.jaeh@mathematik.uni-goettingen.de}
}

\author[Michael Ruzhansky]{Michael Ruzhansky}
\address{
	Michael Ruzhansky:
	\endgraf
	Department of Mathematics: Analysis, Logic
	and Discrete Mathematics
	\endgraf
	Ghent University
	\endgraf
	Belgium
	\endgraf
	and
	\endgraf
	School of Mathematical Sciences
	\endgraf
	Queen Mary University London
	\endgraf
	Mile End, London,  E1 4NS
	\endgraf
	United Kingdom
	\endgraf
	{\it E-mail address} {\rm m.ruzhansky@qmul.ac.uk, \rm michael.ruzhansky@ugent.be}
}

\thanks{
The authors were supported in parts by the FWO Odysseus Project G.0H94.18N: Analysis and Partial Differential Equations, the Leverhulme Grant RPG-2017-151 and by EPSRC Grant
EP/R003025/1.
}

\usepackage{hyperref}


\begin{document}
	
\maketitle

\begin{abstract}
	In this paper we continue the study of non-diagonalisable hyperbolic systems with variable multiplicity started by the authors in \cite{Garetto2018}. In the case of space dependent coefficients, we prove a representation formula for solutions that allows us to derive results of regularity and propagation of singularities.
\end{abstract}
	
\tableofcontents

\section{Introduction}

In this work, we continue the study of non-diagonalisable systems that begun in \cite{Garetto2018} by proving a result on solution representations and propagation of singularities for a hyperbolic system with $x$-dependent principal part. Let us consider \begin{equation} \label{CP_intro} \left\{\begin{aligned}
		& D_t u = A(x,D_x)u + B(t,x,D_x)u + f(t,x), \quad (t,x) \in [0,T] \times \R^n, \\
		& \left. u\right|_{t=0} = u_0, \quad x \in \R^n,
	\end{aligned} \right. \end{equation}  
with the usual notation $D_t=-{\rm i}\partial_t$ and $D_x=-{\rm i}\partial_x$. We assume that $A(x,D_x) = \big( a_{ij}(x,D_x) \big)_{i,j=1}^m$ is an $m\times m$ matrix of pseudo-differential operators of order $1$, i.e., $a_{ij} \in \Psi_{1,0}^1(\R^n))$ and that $B(t,x,D_x) = \big( b_{ij}(t, x,D_x) \big)_{i,j=1}^m$ is an $m\times m$ matrix of pseudo-differential operators of order $0$, i.e., $b_{ij}\in C([0,T],\Psi^0_{1,0}(\R^n))$. We also assume that the matrix $A$ is upper triangular and hyperbolic, i.e., 
	\begin{multline*}
	A(x,D_x) = \Lambda(x,D_x) + N(x,D_x)\\
	={\rm diag}(\lambda_1(x,D_x),\lambda_2(x,D_x),\dots,\lambda_m(x,D_x))+N(x,D_x)
\end{multline*} 
with real eigenvalues $\lambda_1(x,\xi), \lambda_2(x,\xi),\dots, \lambda_m(x,\xi)$ and 
\begin{equation*}
N(x,D_x) = \begin{pmatrix}
	0 & a_{12}(x,D_x) & a_{13}(x,D_x) & \cdots & a_{1m}(x,D_x) \\
	0 & 0 & a_{23}(x,D_x) & \cdots & a_{2m}(x,D_x) \\
	\vdots & \vdots & \vdots & \cdots & \vdots \\
	0 & 0 & 0 & \dots & a_{m-1m}(x,D_x)\\
	0 & 0 & 0 & \dots & 0
\end{pmatrix}.
\end{equation*} 
We recall that the well-posedness of this kind of systems has been proven in anisotropic Sobolev spaces in \cite{Garetto2018} under specific assumptions on the lower order terms.

Propagation of singularities for systems with vanishing iterated Poisson brackets has been studied by several authors as Iwasaki and Morimoto \cite{IwMo84} who studied $3 \times 3$ systems where the twice iterated Poisson  bracket vanishes and Ichinose \cite{Ich82} studied $2 \times 2$ systems under the same condition. In \cite{Roz83}, Rozenblum con\-si\-de\-red smoothly diagonalisable systems with transversally intersecting characteristics, and derived a formula for the propagation of its singularities. Consequently, the transversality condition was removed in
\cite{Kamotski2007}, replaced by a weaker condition of intersection of finite order at points of multiplicity, with propagation of singularities result as well. Here we extend the results of \cite{Kamotski2007} to non-diagonalisable hyperbolic systems with variable multiplicity. The paper is organised as follows. Section 2 collects some important notions of Fourier integral operators and related integral operators relevant to our problem. The main well-posed result and corresponding representation formula is proven in Section 3. Section 4 is devoted to propagation of singularities. The paper ends in Section 5 with an application of our results to higher order hyperbolic equations with multiplicities.

\subsection{Notations and preliminary notions}

For the convenience of the reader we recall here some notations and preliminary notions that we will use throughout the paper.

Let $\mu\in\R$. We recall that $S^\mu_{1,0}(\R^n)$ is the space of symbols of order $\mu$ and type $(1,0)$, i.e., $a=a(x,\xi) \in C^{\infty}(\R^n \times \R^n)$ belongs to $S^\mu_{1,0}(\R^n)$ if there exist constants $C_{\alpha,\beta} >0$ such that \begin{equation*}
		\forall \alpha,\beta \in \N_0^n \,\, : \,\, |\partial_x^\alpha \partial_\xi^\beta a(x,\xi)| \leq C_{\alpha,\beta} \la \xi \ra^{m-|\beta|} \quad \forall (x,\xi) \in \R^n \times \R^n,
	\end{equation*} 
with $\la \xi \ra=(1+|\xi|^2)^{1/2}.$

The set of pseudo-differential operators associated to the symbols in $S^\mu_{1,0}(\R^n)$ is denoted by $\Psi^\mu_{1,0}(\R^n)$. If the symbol has an extra (continuous) dependence on $t\in[0,T]$ we will use the notations $C([0,T], S^\mu)$ and $C([0,T],\Psi^\mu_{1,0})$ for symbols and operators, respectively. For the sake of simplicity we will adopt the abbreviated notations $S^\mu$ and $\Psi^\mu$ for $S^\mu_{1,0}(\R^n)$ and $\Psi^\mu_{1,0}(\R^n)$, respectively, and $CS^\mu$ and $C\Psi^m$ for $C([0,T], S^\mu)$ and $C([0,T],\Psi^\mu_{1,0})$, respectively.

With $I^\mu$ we denote the class of Fourier Integral Operators with amplitude in $S^\mu$, i.e., of operators of the type
\[
I_\varphi(a)(f)(t,x)=\lint_{\R^n} e^{{\rm i}\varphi(t,x,\xi)} a(t,x,\xi) \widehat{f}(\xi) d\xi,
\]
where $\varphi$ is a phase function and $\hat{f}$ is the Fourier transform of $f$. This notation is standard and further details can be found in \cite{Garetto2018} and the references therein, for instance \cite{Hoe71}. In this paper we will use the short expression \emph{integrated Fourier integral operator} to denote an operator of the type
\[
\int_{0}^{t} \lint_{\R^n} e^{{\rm i}\varphi(t,s,x,\xi)} a(t,s,x,\xi) \widehat{f}(\xi,s) d\xi\, ds,
\]
where the Fourier transform of $f=f(y,s)$ is meant with respect to the variable $y$.

By $L_\alpha^p(\R^n)$, we denote the Sobolev space $(I-\Delta)^{-\tfrac{\alpha}{2}} L^p(\R^n)$. As usual, we set $H^s = L_s^2$. By $\|\cdot\|_{L^p_{loc}(\R^n)}$ we denote any localisation of the $L^p(\R^n)$-norm, i.e. the estimate $\|f\|_{L^p_{loc}(\R^n)} \leq C$ means that for all $\chi \in C_0^\infty(\R^n)$, we have the estimate $\|\chi f\|_{L^p(\R^n)} \leq C$, where the constant $C$ may depend on $\chi$. Since we only work in $\R^n$ and no confusion can arise, we drop the indication of $\R^n$ from here on. 

Finally, we recall that the Poisson bracket of two differentiable functions $f=f(x,\xi)$ and $g=g(x,\xi)$ is defined as
\[
\{ f\,,\,g\, \} = \sum_{i=1}^n \frac{\partial f}{\partial \xi_i}\frac{\partial g}{\partial x_i} - \frac{\partial f}{\partial x_i}\frac{\partial g}{\partial \xi_i}.
\]

%We denote by $\hat{f}$ the Fourier transform of $f$ with respect to $x \in \R^n$.

\subsection{Assumptions on the matrices $A(x,D_x)$ and $B(t,x,D_x)$}

In this paper, we make the following assumptions on lower order terms and multiplicities:
 \begin{enumerate}
	\item[(\textbf{H1})] ({\bf Lower order terms}) The entries of the matrix $B(t,x,D_x)=[b_{ij}(t,x,D))]_{i,j=1}^m$ belong to $C([0,T], \Psi^{0})$ and are of decreasing order below the diagonal, i.e.,  

 \begin{equation}
	\label{cond_lo}
	b_{ij}\in C([0,T], \Psi^{j-i})\quad \text{for $i> j$}.
	\end{equation} 
\item[(\textbf{H2})] ({\bf multiplicities}) There exists $M\in\mathbb{N}$ such that if $\lambda_j(x,\xi)=\lambda_k(x,\xi)$ for some $j,k \in \{ 1,\ldots , m \}$ and $\lambda_j(x,\xi)$ and $\lambda_k(x,\xi)$ are not identically equal near $(x,\xi)$ then there exists some $N\le M$ such that \begin{equation*}
	\lambda_j(x,\xi)=\lambda_k(x,\xi) \,\, \Rightarrow \,\, H^N_{\lambda_j}(\lambda_k):=\{\lambda_j,\{\lambda_j,\dots,\{\lambda_j,\lambda_k\}\}\dots\}\neq 0,
	\end{equation*} where the Poisson bracket $\{ \cdot , \cdot \}$ in $H^N_{\lambda_j}$ is iterated $N$ times.
\end{enumerate}

\begin{remark}
	In \cite{Garetto2018}, for $A$ depending on $t$ as well and under the condition (H1) on $B$, we proved that for any $s \in \R$, $u^0_k \in H^{s+k-1}$, $k=1,\ldots,m$, and $f_k \in C([0,T],H^{s+k-1})$, $k=1,\ldots,m$, the Cauchy problem \eqref{CP_intro} has a unique anisotropic Sobolev solution $u$ with components $u_k \in C\left( [0,T], H^{s+k-1} \right)$, $k=1,\ldots,m$. In this paper we do the microlocal analysis of solutions in the case of $A$ depending only on $x$.
\end{remark}

\begin{remark}
	 The condition (H2) was introduced in \cite[p.3]{Kamotski2007}. For $1<p<\infty$ and $\alpha=(n-1)|1/p-1/2|$, it was proved in \cite{Kamotski2007} that when the matrix $A(x,D_x)$ is smoothly microlocally diagonalisable, with smooth eigenspaces and real eigenvalues $\lambda_j$, $j=1,\dots,m$, fulfilling the condition (H2) then for every compactly supported initial data  $u_0\in L^p_\alpha\cap L^2_{comp}$ the Cauchy problem \eqref{CP_intro} has a unique solution $u$ such that $u(t,\cdot)\in L^p_{loc}$ for every $t\in[0,T]$.  
 Previously, a similar result had been proved in $L^2$ by Rozenblum in \cite{Roz83}, in the special case of (H2) with $N=1$.
\end{remark}
We are now ready to state the main result of our paper. This is a representation formula for the solution $u$ which shows how this depends on initial data and right-hand side. This dependence is given in terms of integral operators (of Fourier type) modulo regularising operators of order $N$, i.e. mapping $H^s$ into $H^{s+N}$, for any $s\in\R$. 

\begin{thm} 
\label{thm:SolRep}

Let $n\geq 1$, $m\geq 2$, and let
\[
\left\{\begin{array}{l}
D_t u = A(x,D_x)u + B(t,x,D_x)u + f(t,x), \quad (t,x) \in [0,T] \times \R^n, \\
\left.u \right|_{t=0} = u_0(x), \quad x \in \R^n,
\end{array} \right.
\]
where $A(x,D_x)$ is an upper-triangular matrix of pseudo-diffe\-ren\-tial operators of order $1$ and $B(t,x,D_x)$ is a matrix of pseudo-differential operators of order $0$, continuous with respect to $t$. Let $u_0$ and $f$ have components $u^0_j$ and $f_j$, respectively, with $u^0_j\in H^{s+j-1}(\R^n)$ and $f_j\in C([0,T],H^{s+j-1})$ for $j=1,\dots,m$. Then, under condition (H1) and (H2), we have the following:
\begin{itemize}
\item[(i)] the Cauchy problem above has a unique anisotropic Sobolev solution $u$, i.e.,  $u_j\in C([0,T], H^{s+j-1})$ for $j=1,\dots, m$;
\item[(ii)] for any $N\in\N$, the components $u_j$, $j=1,\ldots,m$, of the solution $u$ are given by
	\begin{equation}
	\label{repr_formula_f}
	u_j(t,x) = \sum_{l=1}^{m} \left( \mc H_{j,l}^{l-j}(t) + R_{j,l}(t) \right) u_l^0 + \left( \mc K_{j,l}^{l-j}(t) + S_{j,l}(t) \right) f_l,
	\end{equation} 
	where $R_{j,l}$, $S_{j,l} \in \mc L(H^s, C([0,T],H^{s+N-l+j}))$ and the operators  $\mc H_{j,l}^{l-j}$, $\mc K_{j,l}^{l-j}$$\in  \mc L(C([0,T], H^s), C([0,T],H^{s-l+j}))$ are integrated Fourier Integral Operators of order $l-j$.
\end{itemize}
\end{thm}
For the convenience of the reader we recall here Theorem 1 in \cite{Garetto2018} which proves already assertion (i) in Theorem \ref{thm:SolRep}.
\begin{thm}[{\bf Theorem 1 in \cite{Garetto2018}}]
	\label{main_theo_mi}
	Let $n\geq 1$, $m\geq 2$, and let
	\begin{equation}
	\label{eq_CP_mi}
	\left\{\begin{array}{l}
	D_t u = A(t,x,D_x)u + B(t,x,D_x)u + f(t,x), \quad (t,x) \in [0,T] \times \R^n, \\
	\left.u \right|_{t=0} = u_0(x), \quad x \in \R^n,
	\end{array} \right.
	\end{equation} where $A(t,x,D_x) = (a_{ij}(t,x,D_x))_{i,j=1}^m$ is an upper-triangular matrix of pseudo-diffe\-ren\-tial operators of order $1$ and $B(t,x,D_x) = (b_{ij}(t,x,D_x))_{i,j=1}^m$  is a matrix of pseudo-differential operators of order $0$,  continuous with respect to $t$. Assume that (H1) holds. Then, given $u^0_j\in H^{s+j-1}(\R^n)$ and $f_j\in C([0,T],H^{s+j-1})$ for $j=1,\dots,m$, the Cauchy problem \eqref{eq_CP_mi} has a unique anisotropic Sobolev solution $u$, i.e.,  $u_j\in C([0,T], H^{s+j-1})$ for $j=1,\dots, m$.
\end{thm}

\section{Auxiliary results} \label{sec:Aux}

This section contains some auxiliary results on Fourier integral operators and related integral operators that we will use throughout the paper. For the convenience of the reader, we begin by recalling some notations introduced in \cite{Garetto2018}.

For each eigenvalue $\lambda_j(x,\xi)$ of $A(x,\xi)$, we will be denoting by $G^0_j\theta$ the solution to \begin{equation*}
\left\{\begin{array}{l}
D_t w = \lambda_j(x,D_x)w + b_{jj}(t,x,D_x)w, \\
w(0,x) = \theta(x),
\end{array}\right.
\end{equation*} and by $G_j g$ the solution to \begin{equation*}
\left\{\begin{array}{l}
D_t w = \lambda_j(x,D_x)w + b_{jj}(t,x,D_x)w + g(t,x), \\
w(0,x) = 0.
\end{array}\right.
\end{equation*} The $b_{jj}$ are the diagonal elements of the lower order term $B(t,x,D_x)$ in \eqref{CP_intro}. The operators $G^0_j$ and $G_j$ can be locally represented by a Fourier integral operator and an integrated  Fourier integral operator, respectively, i.e.,
\begin{equation}\label{eq:Gj0}
G^0_j \theta(t,x) = \lint_{\R^n} e^{{\rm i}\varphi_j(t,x,\xi)} c_j(t,x,\xi) \widehat{\theta}(\xi) d\xi
\end{equation} and \begin{equation} \label{eq:Gj}
G_j g(t,x) = \lint_{0}^{t} \lint_{\R^n} e^{{\rm i}\varphi_j(t,s,x,\xi)} C_j(t,s,x,\xi) \widehat{g}(s,\xi) d\xi ds,
\end{equation} 
with $\varphi_j(t,s,x,\xi)$ solving the eikonal equation
\begin{equation*}
\left\{
\begin{array}{l}
\partial_t \varphi_j = \lambda_j(x,\nabla_x\varphi_j(t,s,x,\xi)), \\
\varphi_j(s,s,x,\xi) = x \cdot \xi,
\end{array}
\right.
\end{equation*} 
and $\varphi_j(t,x,\xi) = \varphi_j(t,0,x,\xi)$. Note that the amplitudes $C_j$ in \eqref{eq:Gj} have asymptotic expansions $\sum\limits_{k=0}^{+\infty} C_{j,-k}$ where the element $C_{j,-k}(s,x,\xi)$ is of order $-k$, $k$ $\in \N$,  and satisfies transport equations with initial data at $t=s$. By construction,  $c_j(t,x,\xi) = C_j(t,0,x,\xi)$. In the above construction of propagators for hyperbolic equations, we have $c_j\in S^0$, so that $G^0_j\in I^{0}.$

Further, to simplify the analysis of the regularising part in \eqref{repr_formula_f} we introduce the notation \begin{equation} \label{def:mcE}
	\mc E_j(t,s)g(s,x) = \lint_{\R^n} e^{{\rm i}\varphi_j(t,s,x,\xi)} C_j(s,x,\xi) \widehat{g}(s,\xi) d\xi,
\end{equation} i.e., the integrated Fourier integral operator $G_j$ can now be written as
\begin{equation}\label{EQ:Ej}
	G_j g(t,x) = \lint_{0}^{t} \mc E_j(t,s)g(s,x) ds.
\end{equation}
 
\subsection{Composition of FIOs and regularising effect}

In this section, we state and prove auxiliary results that are crucial to the proof of the solution representation formula stated in Theorem \ref{thm:SolRep}. In particular, we investigate the mapping properties of compositions and powers of Fourier integral operators and integrated Fourier integral operators as in \eqref{eq:Gj0} and \eqref{eq:Gj}. This will be useful when analysing the regularising part of our representation formula \eqref{repr_formula_f}.

\subsubsection{Integrated Fourier integral operators}

The composition of integrated Fourier integral operators like in \eqref{eq:Gj} was studied in \cite{Kamotski2007}. Their result, that we recall in the sequel, is crucial for our proof and is a generalisation of a previous result by Rozenblum in \cite{Roz83}.
  
Let $t_1, t_2, \dots, t_l\in[0,T]$, $\overline{t} = (t_1,t_2,\ldots,t_l)$ and let $H(\overline{t})$ be the operator \begin{equation}\label{def:HFIO}
H(\overline{t}) = e^{i\lambda_{j_1}t_1}e^{i\lambda_{j_2}(t_2-t_1)} \cdot \ldots \cdot e^{i\lambda_{j_l}(t_{l}-t_{l-1})} e^{-i\lambda_{j_l}t_{l}},
\end{equation} where $\lambda_i$ are pseudo-differential operators of order $1$.   

%\textcolor{blue}{???? Note that the exponentials in \eqref{def:HFIO} correspond to the $\mathcal E_j$ in \eqref{def:mcE} where the pseudo-differential terms coming from the nilpotent part of $A$ and the lower order term $B$, have been 'moved out' by means of Egorov's theorem; that is the $B(\overline{t})$ in Theorem \ref{thm:RegFIO}.}

%\begin{remark}
%	By Egorov's theorem, we can consider $G_i \circ G_j$ instead of $G_i \circ A \circ G_j B$ since \begin{equation*}
%	G_i \circ A \circ G_j B = G_i \circ G_j \circ G_j^{-1} \circ A \circ G_j \circ B = G_i G_j C.
%	\end{equation*}
%\end{remark}

By \cite{Hoe71}, $H(\overline{t})$ is a (parameter dependent) Fourier integral operator and its canonical relation $\Lambda^{\overline{t}} \subseteq T^\ast \R^n \times T^\ast \R^n$ is given by \begin{equation*}
\Lambda^{\overline{t}} = \left\{ (x,p,y,\xi) \,:\,\, (x,p) = \Psi^{\overline{t}}(y,\xi) \right\},
\end{equation*} where \begin{equation*}
\Psi^{\overline{t}} = \Phi_{j_1}^{t_1} \circ \dots \circ \Phi_{j_l}^{t_l-t_{l-1}} \circ \Phi_{j_{l+1}}^{-t_l}
\end{equation*} and the $\Phi_j^t$ are the transformations corresponding to a shift by $t$ along the trajectories of the Hamiltonian flow defined by the $\lambda_j$.

\begin{thm}[Thm 2.1 in \cite{Kamotski2007}] \label{thm:RegFIO} With the above notation, assume that not all $\lambda_i$s are identical to each other and let (H2) be satisfied for the $\lambda_j$ in \eqref{def:HFIO}. Further, suppose that $D(\overline{t}) \in \Psi^0$. Then, the operator \begin{equation*}
	Q_l = \lint_{0}^t \lint_{0}^{t_1} \ldots \lint_{0}^{t_{l-1}} D(\overline{t})H(\overline{t}) \, dt_{l} \ldots dt_1
	\end{equation*} belongs to $\mc L(H^s,H^{s+N(l)})$, where $N(l) \to +\infty$ as $l \to +\infty$.
\end{thm}

\begin{remark} \label{rem:RKR}
	If the global estimate in the definition of the symbols classes $S^m_{1,0}$ in \cite{Garetto2018} is replaced with an estimate that holds locally on every compact set, then the the conclusions of Theorem \ref{thm:RegFIO} hold true if $\mc L(H^s,H^{s+N(l)})$ is replaced by $\mc L(H^s_{comp}, H^{s+N(l)}_{loc})$.
\end{remark}

Theorem \ref{thm:RegFIO} allows us to investigate the composition $G_iG_j$.

\subsubsection{The composition $G_iG_j$}

Let $G_i$ and $G_j$ be two operators as in \eqref{eq:Gj}. Then, we can write \begin{equation*}
	G_i G_j u = \lint_{0}^t \lint_0^{t_1} \mc E_i(t,t_1) \mc E_j(t_1,t_2)u \, dt_2 dt_1,
\end{equation*} 
with $\mc E_i, \mc E_j$ as in \eqref{EQ:Ej}.
If we now iterate this $k$ times, we obtain \begin{equation*}
	(G_i G_j)^k  = \lint_{0}^t \lint_0^{t_1} \dots \lint_0^{t_{2k-1}} \mc E_{i}(t,t_1)\mc E_{j}(t_1,t_2)\mc \cdot \ldots \cdot \mc E_i(t_{2k-3},t_{2k-2}) \mc E_j(t_ {2k-2},t_{2k-1}) \,d\overline{t},
\end{equation*} where $\overline{t} = dt_1 \ldots dt_{2k-1}$. This is an operator of the same type as $Q_l$ above so we can apply Theorem \ref{thm:RegFIO} and obtain the following corollary.

\begin{cor}
\label{cor_regularising}
Let conditions (H1) and (H2) be satisfied. Let $j \in \{ 1,\ldots, m \}$, $s \in \R$. Then, for all $N \in \N$, $k \in \{ 1,\ldots ,m \}$, there exists $M \in \N$ such that 
\begin{equation*}
	\left( \prod_{i=1}^k G_{\sigma(i)} \right)^M
	 \in \mc L\left(C([0,T],H^s),C([0,T],H^{s+N})\right),
	\end{equation*} where $\sigma$ is an element of the symmetric group over $\{1,\ldots,m\}$.
\end{cor}

\begin{remark} \label{ren:ExtCorRegFIO}
	The same conclusion as in Corollary  \ref{cor_regularising} holds true if we have a product that contains a collection $G_j^0$'s as long as there is at least one integrated version $G_j$ present.
\end{remark}

%\begin{remark}
%	Note that in \cite{Garetto2018}, we defined $S^m$ by global estimates on $\R^n \times \R^n$. If  one replaces that definition with a locally over every compact sets version, Theorem \ref{thm:SolRep} still holds true with the spaces $\mc L(H^s, C([0,T],H^{s+N-l+j}))$ and $\mc L(C([0,T], H^s), C([0,T],H^{s-l+j}))$ replaced by $\mc L(H^s_{comp}, C([0,T],H^{s+N-l+j}_{loc}))$ and\\ 
%$\mc L(C([0,T], H^s_{comp}), C([0,T],H^{s-l+j}_{loc}))$, respectively.
%\end{remark}

\section{Solution representations}

This section is devoted to the proof of Theorem \ref{thm:SolRep}. For the sake of simplicity and for the advantage of the reader we give first a detailed explanatory proof for $2\times 2$ systems and we then pass to consider the $m\times m$ case. We adopt the notations introduced in Section \ref{sec:Aux}.   
\subsection{The $2\times 2$ case}

Let us consider the system \begin{eqnarray} \label{eq:CP2x2}
\left\{ \begin{array}{l}
D_t u = A(x,D_x)u + B(t,x,D_x)u + f(t,x), \quad (t,x) \in [0,T] \times \R^n,\\
\left.u\right|_{t=0} = u_0, \quad x \in \R^n,
\end{array}
\right.
\end{eqnarray} where $u_0(x) = \big(u_1^0(x),u_2^0(x) \big)^T$, $f(t,x) = \big( f_1(t,x) ,f_2(t,x) \big)^T$, and with the operators $A(x,D_x)$ and $B(t,x,D_x)$  given by \begin{equation} \label{eq:A2x2}
A(x,D_x) = \begin{pmatrix}
\lambda_1(x,D_x) & a_{12}(x,D_x) \\
0 & \lambda_2(x,D_x) \\
\end{pmatrix} \end{equation} 
and 
\begin{equation*}
B(t,x,D_x) = \begin{pmatrix}
b_{11}(t,x,D_x) & b_{12}(t,x,D_x) \\
b_{21}(t,x,D_x) & b_{22}(t,x,D_x) \\
\end{pmatrix}.
\end{equation*} We suppose that all entries of $A(x,D_x)$ belong to $\Psi_{1,0}^1$ and all entries of $B(t,x,D_x)$ belong to $ C \Psi_{1,0}^0$.

As detailed in Subsection 2.2 in \cite{Garetto2018}, we obtain
the equations \begin{eqnarray*}
	u_1 &=& G_1^0(u_1^0) + G_1(f_1) + G_1(a_{12}u_2) + G_1(b_{12}u_2) \\
	&=& G_1^0(u_1^0) + G_1(f_1) + G_1((a_{12}+b_{12})u_2), \\
	u_2 &=& G_2^0(u_2^0) + G_2(f_2) + G_2(b_{21}u_1),
\end{eqnarray*} and with that \begin{equation} \label{eq:Solu1}
\begin{aligned}
u_1 &= G_1^0(u_1^0) + G_1((a_{12}+b_{12})G_2^0(u_2^0)) +  G_1(f_1) + G_1((a_{12}+b_{12})G_2(f_2))  \\
& \qquad  + G_1((a_{12}+b_{12})G_2(b_{21}u_1)), \\
u_2 &= G_2^0(u_2^0) + G_2(b_{21}G_1^0(u_0^1)) + G_2(b_{21}G_1(f_1)) + G_2(f_2) \\
& \qquad + G_2(b_{21}G_1((a_{12}+b_{12})u_2)).
\end{aligned}
\end{equation}

We note that the operators  $G_1 \circ (a_{12}+b_{12}) \circ G_2 \circ b_{21}$ and $G_2 \circ b_{21} \circ G_1 \circ (a_{12}+b_{12})$ are of order $0$ under the assumption (H1) made on the lower order terms; here in particular $b_{21} \in \Psi^{-1}$. From \eqref{eq:Solu1}, we have \begin{equation*}
\begin{aligned}
&u_1 - G_1((a_{12}+b_{12})G_2(b_{21}u_1)) \\
&\qquad \qquad = G_1^0(u_1^0) + G_1((a_{12}+b_{12})G_2^0(u_2^0))  + G_1(f_1) + G_1((a_{12}+b_{12})G_2(f_2))
\end{aligned}
\end{equation*} and \begin{equation} \label{eq:u2_resolved}
\begin{aligned}
& u_2 - G_2(b_{21}G_1((a_{12}+b_{12})u_2)) \\
& \qquad \qquad = G_2^0(u_2^0) + G_2(b_{21}G_1^0(u_0^1)) + G_2(b_{21}G_1(f_1)) + G_2(f_2).
\end{aligned}
\end{equation}

\subsection{Inversion of the operator $L_1$}
Adopting the notations introduced in \cite{Garetto2018} we introduce the operator 
\begin{equation*}
L_1 := \textbf{I} - G_1 \circ (a_{12}+b_{12}) \circ G_2 \circ b_{21} = \textbf{I}- \mc G^0_1
\end{equation*} 
Note that $\mc G^0_1=G_1 \circ (a_{12}+b_{12})\circ G_2 \circ b_{21}$ is of order $0$ and and from the Sobolev mapping properties of Fourier Integral Operators (see Lemma 1 in \cite{Garetto2018}) the norm of this operator can be estimated by a constant times the length of the time interval $[0,T]$. So it can be made as small as wanted by a suitable choice of $T$. It follows that $L_1$ is invertible for $T$ small enough and its inverse can be written as sum of a Neumann series. More precisely, under the assumptions (H1) and (H2) from Corollary \ref{cor_regularising} we get that for every $N\in\N$ the operator $L^{-1}$ can be written as a finite sum of powers of the operator $\mc G^0_1$ modulo some regularising operator mapping $C([0,T],H^s)$ into $C([0,T],H^{s+N}))$, i.e., for every $N \in \N$, there exists $M \in N$ such that 
\begin{eqnarray*}
	L_1^{-1}  &=& \sum\limits_{k=0}^{+\infty} \left(\mc G^0_1 \right)^k \\
	&=& \sum\limits_{k=0}^{M} \left(\mc G^0_1 \right)^k \mod \, \mathcal L(C([0,T],H^s),C([0,T],H^{s+N})).
%	&=& \underbrace{\mc H}_{\text{order 0}} + \underbrace{\sum_{k=N+1}^{+\infty} \left(\mc G^0_1 \right)^k}_{_{\text{order -N}}}
\end{eqnarray*}

It is important to remark here that the estimates needed to ensure the small norm of the operator $\mc G_1^0$, do not depend on the initial data and therefore one can repeat the same argument covering the original interval $[0,T]$.

%Question: {\color{blue} How do we see that the remainder has the desired regularising properties? The hope is that it behaves like the operator $Q_l$ in \cite[Thm. 2.1]{Kamotski2007} since the many time integrals may be used as in \cite[p. 11]{Kamotski2007}. The change of variables $\bar{t} = \zeta |\xi|^{-1}$ transforms $Q_l$ into a FIO with integration over a cone. See also Section 4 in \cite{Kamotski2007}.}

%\begin{remark}
%	That we can actually do this remains to be justified. One needs to follow the analysis of \cite[Thm. 2.1]{Kamotski2007}; see also \cite[p. 11 and pp. 19]{Kamotski2007}. In the case of just $x$-dependent coefficients and $N \equiv 0$, we should be able to apply the analysis directly. In the presence of off-diagonal terms, we will have to justify that the analysis still goes through.
%	In the $t$-dependent case, the condition C in \cite[p. 3]{Kamotski2007} has to be modified accordingly.
%\end{remark}
\subsection{Representation formulas}
We now apply the operator $L_1^{-1}$ as written above to both sides of the equality
\[
L_1u_1 = G_1^0(u_1^0) + G_1((a_{12}+b_{12})G_2^0(u_2^0))  + G_1(f_1) + G_1((a_{12}+b_{12})G_2(f_2)).
\]
We obtain the following representation for $u_1$, where
\[
R_1=L_1^{-1} - \sum\limits_{k=0}^{M} \left(\mc G^0_1 \right)^k
\] 
is a regularising operator, i.e., $R_1\in\mathcal L(C([0,T],H^s),C([0,T],H^{s+N}))$:
 \begin{eqnarray*}
	u_1 &=& \underbrace{\sum\limits_{k=0}^{M} \left(\mc G^0_1 \right)^k G_1^0(u_1^0)}_{= \mc H_{1,1}^{1-1}u_1^0} + \underbrace{\sum\limits_{k=0}^{M} \left(\mc G^0_1 \right)^k G_1((a_{12}+b_{12})G_2^0(u_2^0))}_{= \mc H^{2-1}_{1,2} u_2^0} \\
	&& \quad + \underbrace{\sum\limits_{k=0}^{M} \left(\mc G^0_1 \right)^k G_1(f_1)}_{= \mc K_{1,1}^{1-1}f_1} + \underbrace{\sum\limits_{k=0}^{M} \left(\mc G^0_1 \right)^k G_1((a_{12}+b_{12})G_2(f_2))}_{= \mc K^{2-1}_{1,2} f_2} \\
	&& \quad + R_1 G_1^0(u_1^0) + R_1 G_1((a_{12}+b_{12})G_2^0(u_2^0)) \\
	&& \quad + R_1 G_1(f_1) + R_1 G_1((a_{12}+b_{12})G_2(f_2)).
\end{eqnarray*} 
Denoting $R_1 G_1^0$ and $R_1 G_1((a_{12}+b_{12})G_2^0$ by $R_{1,1}$ and $R_{1,2}$ respectively, and  $R_1 G_1$ and $R_1 G_1((a_{12}+b_{12})G_2$ by $S_{1,1}$ and $S_{1,2}$, respectively, we have that
\[
u_1 = \sum_{l=1}^{2} \left( \mc H_{1,l}^{l-j}(t) + R_{1,l}(t) \right) u_l^0 + \left( \mc K_{1,l}^{l-j}(t) + S_{1,l}(t) \right) f_l,
\]
where 
\begin{itemize}
\item the operators $H_{1,l}^{l-1}$ and $K_{1,l}^{l-1}$ are of order $l-1$ and therefore map $C([0,T], H^s)$ into $C([0,T],H^{s-l+1})$,
\item $R_{1,1}$ and $S_{1,1}$ map $H^s$ into $C([0,T],H^{s+N})$,
\item $R_{1,2}$ and $S_{1,2}$ map $H^s$ into $C([0,T],H^{s+N-1})$.
\end{itemize}
This means that $\mc H_{1,l}^{1-l}$, $\mc K_{1,l}^{1-l}$$\in  \mc L(C([0,T], H^s), C([0,T],H^{s-l+1}))$ and $R_{1,l}$, $S_{1,l} \in \mc L(H^s, C([0,T],H^{s+N-l+1}))$. The same argument is true for $u_2$. We have in this way obtained the representation formula stated in Theorem \ref{thm:SolRep}.

\subsection{The $m \times m$ case}
We are now ready to prove Theorem \ref{thm:SolRep}. 
\begin{proof}[Proof of Theorem \ref{thm:SolRep}]
Throughout this proof we refer to the proof of Theorem 1 in \cite{Garetto2018}, i.e. Theorem \ref{main_theo_mi} in this paper. Theorem \ref{main_theo_mi} proves the well-posedness of this Cauchy problem in anisotropic Sobolev spaces. It remains to prove the representation formula for the components of the solution $u$. We begin by observing that under our hypotheses we can write
\[
u_i=U_i^0+\sum_{j<i} G^{j-i}_{i,j}(u_j)+\sum_{i<j\le m} G^1_{i,j}(u_j),
\]
for $i=1,\dots,m$, where $G^{j-i}_{i,j}$ and $G^1_{i,j}$ are operators of order $j-i$ and $1$, respectively and 
\[
U_i^0 = G_i^0 u_j^0 + G_i(f_i).
\]
We begin by substituting  
\[
u_m=U_m^0+\sum_{j<m} G^{j-m}_{m,j}(u_j),
\]
into
\[
u_{m-1}=U_{m-1}^0+\sum_{j<m-1} G^{j-m+1}_{m-1,j}(u_j)+G^1_{m-1,m}(u_m).
\]
We get
\[
\begin{split}
u_{m-1}&=U_{m-1}^0+\sum_{j<m-1} G^{j-m+1}_{m-1,j}(u_j)+G^1_{m-1,m}U_m^0+\sum_{j<m} G^1_{m-1,m}G^{j-m}_{m,j}(u_j)\\
&=(U_{m-1}^0+G^1_{m-1,m}U_m^0)+\sum_{j<m-1} (G^{j-m+1}_{m-1,j}(u_j)+G^1_{m-1,m}G^{j-m}_{m,j}(u_j))\\
&\qquad + G^1_{m-1,m}G^{-1}_{m,m-1}u_{m-1}.
\end{split}
\]
Since all the operators above are of order $\le 0$ we conclude that the operator
\[
L_{m-1}=I- G^1_{m-1,m}G^{-1}_{m,m-1}:=I-\mathcal{G}^0_{m-1}
\]
is invertible on a sufficiently small interval $[0,T]$ and, therefore, 
\begin{equation}
\label{u_m-1}
\begin{aligned}
& u_{m-1}- G^1_{m-1,m}G^{-1}_{m,m-1}u_{m-1} \\
& \qquad =(U_{m-1}^0+G^1_{m-1,m}U_m^0) \\
& \qquad \qquad +\sum_{j<m-1} (G^{j-m+1}_{m-1,j}(u_j)+G^1_{m-1,m}G^{j-m}_{m,j}(u_j)),
\end{aligned}
\end{equation}
yields
\begin{equation}
\label{L_m-1}
u_{m-1}=L^{-1}_{m-1}\widetilde{U}^0_{m-1}+L^{-1}_{m-1}\sum_{j<m-1}\widetilde{G}^{j-m+1}_{m-1}u_j,
\end{equation}
with $\widetilde{U}^0_{m-1}$ and $\widetilde{G}^{j-m+1}_{m-1}$ defined by the right-hand side of \eqref{u_m-1}. In particular,
\begin{equation}
\label{U0m-1}
\widetilde{U}^0_{m-1}=U_{m-1}^0+D^1_{m-1,m}U_m^0,
\end{equation}
where $G^1_{m-1,m}=D^1_{m-1,m}$ is an integrated Fourier integral operator with symbol of order $1$. Note that we choose the notation $D^1_{m-1,m}$ in order to have simpler notations for the compositions of operators in the computations below. Substituting $u_{m}$ and $u_{m-1}$ into $u_{m-2}$ and making use of \eqref{L_m-1} we find a similar formula to \eqref{L_m-1} for $u_{m-2}$ (see (24) in \cite{Garetto2018}) with $\widetilde{U}^0_{m-2}$ defined as follows:
\begin{equation}
\label{indata_m-2}
\begin{aligned}
& \widetilde{U}^0_{m-2}=U_{m-2}^0+G^1_{m-2,m-1}L^{-1}_{m-1}\widetilde{U}^0_{m-1}\\
& \qquad \qquad +G^1_{m-2,m}U^0_m+G^1_{m-2,m}G^{-1}_{m,m-1}L^{-1}_{m-1}\widetilde{U}^0_{m-1}.
\end{aligned} 
\end{equation}
Hence, by implementing \eqref{U0m-1} in \eqref{indata_m-2} we have
\[
\begin{aligned}
& \widetilde{U}^0_{m-2}=U_{m-2}^0+G^1_{m-2,m-1}L^{-1}_{m-1}(U_{m-1}^0+D^1_{m-1,m}U_m^0)\\
& \qquad \qquad +G^1_{m-2,m}U^0_m+G^1_{m-2,m}G^{-1}_{m,m-1}L^{-1}_{m-1}(U_{m-1}^0+D^1_{m-1,m}U_m^0).
\end{aligned} 
\]
By collecting the terms $U^0_{m-1}$ and $U^0_{m}$ we conclude that $\widetilde{U}^0_{m-2}$ can be written as
\[
\widetilde{U}^0_{m-2}=U^0_{m-2}+D^1_{m-2,m-1}U^0_{m-1}+D^2_{m-2,m}U^0_m,
\]
where the operators $D^1_{m-2,m-1}$ and $D^2_{m-2,m}$ are of order $1$ and $2$, respectively. By iterating the same argument we prove that for every $j=1,\dots,m-1$,
\[
\widetilde{U}^0_{j}=U^0_j+\sum_{k>j}D^{k-j}_{j,k}U^0_k,
\]
where $k-j$ is the order of the operator $D^{k-j}_{j,k}$. For a precise construction of the operators $D^{k-j}_{j,k}$ we refer the reader to the proof of Theorem 1 in \cite{Garetto2018}. Since 
\[
u_1= \widetilde{U}^0_{1}+\mathcal{G}^0_{1}u_1,
\]
where the operator $I-\mathcal{G}^0_{1}$ is invertible with inverse $L_1^{-1}$ (see \cite{Garetto2018}) we conclude that
\[
u_1=L_1^{-1}\biggl(U^0_1+\sum_{k>1}D^{k-1}_{1,k}U^0_k\biggr)=L_1^{-1}\biggl(G_1^0 u_1^0 + G_1(f_1)+\sum_{k>1}D^{k-1}_{1,k}(G_k^0 u_k^0 + G_k(f_k))\biggr).
\]
We now argue as in the case $2\times 2$ and we apply Corollary \ref{cor_regularising} which holds thanks to the hypotheses (H1) and (H2). We obtain that for every $N$ there exists $M$ such that $L_1^{-1}=\sum_{j=0}^M(\mathcal G_1^0)^j$ modulo some regularising operator $R_1$ mapping $C([0,T],H^s)$ into $C([0,T],H^{s+N})$ and therefore
\begin{multline*}
u_1=\sum_{j=0}^M(\mathcal G_1^0)^j\biggl(G_1^0 u_1^0 + G_1(f_1)+\sum_{k>1}D^{k-1}_{1,k}(G_k^0 u_k^0 + G_k(f_k))\biggr)\\
+R_1\biggl(G_1^0 u_1^0 + G_1(f_1)+\sum_{k>1}D^{k-1}_{1,k}(G_k^0 u_k^0 + G_k(f_k))\biggr).
\end{multline*}
By collecting all the terms with $u_k^0$ and all the terms with $f_k$, $k=1,\dots,m$ we see that 
\[
u_1 = \sum_{l=1}^{m} \left( \mc H_{1,l}^{l-1}(t) + R_{1,l}(t) \right) u_l^0 + \left( \mc K_{1,l}^{l-1}(t) + S_{1,l}(t) \right) f_l,
\]
where $\mc H_{1,l}^{l-1}$ and $\mc K_{1,l}^{l-1}$ have order $l-1$ and $R_{1,l}$ and $S_{1,l}$ are regularising.
This is due to the fact $(\mathcal G_1^0)^j$ is an operator of order $0$ as well as $G_k^0$ and $G_k$, and $D^{k-1}_{1,k}$ is an operator of order $k-1$. The regularising operator $R_1$ generates $R_{1,l}$ and $S_{1,l}$. These last two operators map $C([0,T], H^s)$ into $C([0,T],H^{s+N-l+1})$. We have therefore proven the second assertion of this theorem for $j=1$. Following the proof of Theorem 1 in \cite{Garetto2018} we now have that 
\[
u_2=\widetilde{U}^0_{2}+\mathcal{G}^0_{2}u_2+\widetilde{G}^{-1}_{2}u_1,
\]
where the operator $\mathcal{G}^0_2$ is of zero order and its definition involves invertible operators $L_{m-1}, L_{m-2}, \ldots ,L_{2}$ and $\widetilde{G}^{-1}_{2}$ is of order $-1$. Hence, by inverting the operator $L_2=I-\mathcal{G}^0_2$ on a sufficiently small interval $[0,T]$ we have
\[
u_2=L^{-1}_2\widetilde{U}^0_{2}+L^{-1}_2\widetilde{G}^{-1}_{2}u_1.
\]
By definition of $\widetilde{U}^0_{2}$ and by the representation formula for $u_1$ obtained above we can write
\begin{multline*}
u_2=L_2^{-1}\biggl(G_2^0 u_2^0 + G_2(f_2)+\sum_{k>2}D^{k-2}_{2,k}(G_k^0 u_k^0 + G_k(f_k))\biggr)\\
+L_2^{-1}\widetilde{G}^{-1}_{2}\biggl(\sum_{l=1}^{m} \left( \mc H_{1,l}^{l-1}(t) + R_{1,l}(t) \right) u_l^0 + \left( \mc K_{1,l}^{l-1}(t) + S_{1,l}(t) \right) f_l\biggr).
\end{multline*}
Note that the operators above are of order $l-2$. By arguing as for $u_1$ and by writing $L^{-1}_2$ as a finite number of powers of $\widetilde{G}^{-1}_{2}$ plus a regularising operator we arrive at the formula
\[
u_2 = \sum_{l=1}^{m} \left( \mc H_{2,l}^{l-2}(t) + R_{2,l}(t) \right) u_l^0 + \left( \mc K_{2,l}^{l-2}(t) + S_{2,l}(t) \right) f_l,
\] 
with the desired order and regularising properties. We conclude the proof by iterating the same scheme. From formula (28) in \cite{Garetto2018} we obtain for $j>2$ the following expression for $u_j$:
\begin{multline*}
u_j=L_j^{-1} \biggl(U^0_j+\sum_{k>j}D^{k-j}_{j,k}U^0_k\biggr)+L_j^{-1}\biggl(\sum_{k<j}\widetilde{G}_j^{k-j}u_k\biggr)\\
=L_j^{-1}\biggl(G_j^0 u_j^0 + G_j(f_j)+\sum_{k>j}D^{k-j}_{j,k}(G_k^0 u_k^0 + G_k(f_k))\biggr)\\
+L_j^{-1}\biggl(\sum_{k<j}\widetilde{G}_j^{k-j}\biggl( \sum_{l=1}^{m} \left( \mc H_{k,l}^{l-k}(t) + R_{k,l}(t) \right) u_l^0 + \left( \mc K_{k,l}^{l-k}(t) + S_{k,l}(t) \right) f_l\biggr)\biggr),
\end{multline*}
where the operator involved are of order $l-j$. Writing $L_j^{-1}$ by Neumann series we can conclude that
\[
u_j = \sum_{l=1}^{m} \left( \mc H_{j,l}^{l-j}(t) + R_{j,l}(t) \right) u_l^0 + \left( \mc K_{j,l}^{l-j}(t) + S_{j,l}(t) \right) f_l,
\]
where $\mc H_{j,l}^{l-j}$ and $\mc K_{j,l}^{l-j}$ have order $l-j$ and $R_{j,l}$ and $S_{j,l}$ map $C([0,T], H^s)$ into $C([0,T],H^{s+N-l+j})$.
\end{proof} 

\begin{remark}
	Note that in \cite{Garetto2018}, we defined $S^m$ by global estimates on $\R^n \times \R^n$. If  one replaces that definition with a locally over every compact sets version then Theorem \ref{thm:SolRep} still holds true with the spaces $\mc L(H^s, C([0,T],H^{s+N}))$ and $\mc L(C([0,T], H^s), C([0,T],H^{s+N-l+j}))$  replaced by the spaces $\mc L(H^s_{comp}, C([0,T],H^{s+N}_{loc}))$ and 
$\mc L(C([0,T], H^s_{comp}, C([0,T],H^{s+N-l+j}_{loc}))$, respectively.
\end{remark}

%In \cite{Seeger91}, it was shown that a Fourier integral operator of order zero satisfying the local graph condition, is locally bounded from $(L_\alpha^p)_{comp}$ to $L^p_{loc}$ for $p \in (1,+\infty)$ and $\alpha=(n-1)\big|\tfrac{1}{p}-\tfrac{1}{2}\big|$. As stated in \cite{Kamotski2007}, one readily deduces that if \eqref{CP_intro} is strictly hyperbolic, there 

\subsection{Regularity results}

We conclude this section with some regularity results in $L^p$ and H\"older spaces. These are obtained by arguing as in \cite{Kamotski2007} Theorem 2.2 and Theorem 3.1.

\begin{thm}
	Let $1<p<\infty$ and $\alpha=(n-1)\big|\tfrac{1}{p}-\tfrac{1}{2}\big|$. Let $A(x,D_x)$ be an $m \times m$ upper-triangular matrix of pseudo-differential operators of order $1$ and suppose that the eigenvalues $\lambda_i(x,\xi) \in S^1$ of $A(x,\xi)$ are real and satisfy (H2). Assume further, that $B(t,x,D_x)$ is an $m \times m$ matrix of pseudo-differential operators of order $0$ satisfying (H1). Then, for any compactly supported $u_0 \in L_\alpha^p \cap L^2_{comp}$, the solution $u=u(t,x)$ of the Cauchy problem \eqref{CP_intro} satisfies $u(t,\cdot) \in L^p_{loc}$, for all $t \in [0,T]$. Moreover, there is a positive constant $C_T$ such that \begin{equation*}
		\sup_{t \in [0,T]} \|u(t,\cdot)\|_{L^p_{loc}} \leq C_T \|u_0\|_{L_\alpha^p}.
	\end{equation*}
\end{thm}

Local estimates can be obtained in other spaces as well, for $s\in\R$ and $\alpha$ as above. In detail, assuming $u_0$ below is compactly supported, we have \begin{itemize}
	\item $u_0 \in L_{s+\alpha}^p$ implies $u(t,\cdot) \in L_{s}^p$;
	\item $u_0 \in C^{s+\frac{n-1}{2}}$ implies $u(t,\cdot) \in C^s$;
	\item for $1<p\leq q \leq 2$, $u_0 \in L^p_{s-\tfrac{1}{q}+\tfrac{n}{p}-\tfrac{n-1}{2}}$ implies $u(t,\cdot) \in L^q_s$.
\end{itemize}

\section{Propagation of singularities}
We now want to analyse the solution $u$ under a microlocal point of view. In particular we want to see how its wavefront set is related to the wavefront set of the initial data. Thanks to the assumptions (H1) and (H2) and the representation formula in Theorem \ref{thm:SolRep} we are able to extend the result of propagation of singularities in  \cite{Kamotski2007} to systems with not diagonalisable principal part (in upper-triangular form). For the sake of the reader we recall below some basic notions of microlocal analysis which can be found in \cite{Hoe90} and \cite{Hoe71}.
\begin{definition}[Def. 8.1.2 in \cite{Hoe90}, Def. 2.5.2 in \cite{Hoe71}]
	Let $v \in \mc D'(\R^n)$. The wave front set $\WF(v) \subseteq T^\ast(\R^n)\setminus\{0\}:= \R^n \times \R^n\setminus\{0\}$ is defined via its complement as follows: $(x_0,\xi_0)$ belongs to $(\WF(u))^c$ if and only if there exists a $\chi \in C_0^\infty(\R^n)$ with $\chi(x_0) \neq 0$ and a conic neighbourhood $\Gamma$ of $\xi_0$ such that for every $N \in \N$ there exists a positive constant $C_{N,\chi}$ such that \begin{equation*}
	  |\mathcal{F}(\chi u)(\xi)| \leq C_{n,\chi} \la \xi \ra^{-N},
	\end{equation*}
	for all $\xi \in \Gamma$. 
\end{definition}

Let us now discuss the propagation of singularities 
for operators $Q_l$ from
Theorem \ref{thm:RegFIO}, given by 
\begin{equation}\label{EQ:Ql}
	Q_l = \lint_{0}^t \lint_{0}^{t_1} \ldots \lint_{0}^{t_{l-1}} D(\overline{t})H(\overline{t}) \, dt_{l} \ldots dt_1.
\end{equation}
Similar analysis was done in \cite{Kamotski2007}.
It is clear that singularities propagate
along broken Hamiltonian flows. Let
$$J=\{j_1,\ldots,j_{l+1}\}, \; 1\leq j_k\leq m,\; j_k\not=j_{k+1}.$$
We recall from the definition of $H(\overline{t})$ that its
canonical relation $\Lambda^{\overline{t}} \subseteq T^\ast \R^n \times T^\ast \R^n$ is given by \begin{equation*}
\Lambda^{\overline{t}} = \left\{ (x,p,y,\xi) \,:\,\, (x,p) = \Psi^{\overline{t}}(y,\xi) \right\},
\end{equation*} where \begin{equation*}
\Psi^{\overline{t}} = \Phi_{j_1}^{t_1} \circ \dots \circ \Phi_{j_l}^{t_l-t_{l-1}} \circ \Phi_{j_{l+1}}^{-t_l}
\end{equation*} and the $\Phi_j^t$ are the transformations corresponding to a shift by $t$ along the trajectories of the Hamiltonian flow defined by the $\lambda_j$.
Also, recall that $t_1, t_2, \dots, t_l\in[0,T]$, $\overline{t} = (t_1,t_2,\ldots,t_l)$, and $H(\overline{t})$ is the operator \begin{equation}\label{def:HFIO2}
H(\overline{t}) = e^{i\lambda_{j_1}t_1}e^{i\lambda_{j_2}(t_2-t_1)} \cdot \ldots \cdot e^{i\lambda_{j_l}(t_{l}-t_{l-1})} e^{-i\lambda_{j_l}t_{l}},
\end{equation} 
where $\lambda_i$ are pseudo-differential operators of order $1$. 

Let $\Phi_J(t,x,\xi)$
be the corresponding broken Hamiltonian flow. It means that points
follow bicharacteristics of $\lambda_{j_1}$ until meeting the characteristic
of $\lambda_{j_2}$, and then continue along the bicharacteristic of $\lambda_{j_2}$,
etc.
In this procedure the singularities may accumulate if wave front sets for different
broken trajectories project to the same point of $X$.
We can rewrite \eqref{EQ:Ql} as
$$Q_l=\int_\Delta I(\bar{t}) d\bar{t},$$
where $\bar{t}=(t_1,\ldots,t_l)$ ranges over the simplex
$$\Delta=\{0\leq t_l\leq t_{l-1}\leq \ldots\leq t_1\leq t\}$$
in $\Rone^l$ and
$$I(\bar{t})=Z(t_1)\circ\ldots\circ Z(t_l),$$
with $Z(t_j)$ found from \eqref{def:HFIO2}.
It is then possible to treat it as
a standard Fourier integral operator
with the change of variables $\bar{t}=\zeta|\xi|^{-1}.$
Let $K$ be a cone in $\Rone^N=\Rone^{n+l}.$
Let
$$Iu(x)=\int_K\int_Y e^{i\varphi(x,y,\theta)} a(x,y,\theta) u(y) dy d\theta$$
be a Fourier integral operator with integration over the cone $K$ with respect
to $\theta.$
Let $K_j$ be $K$ or a face of $K$. Let $\varphi_j(x,y,\theta_j)=\varphi|_{K_j},
\theta_j\in K_j.$ Let $\Lambda_j\subset T^*X\times T^*X$ be a
Lagrangian manifold with boundary:
$$\Lambda_j=\left\{\left(x,\frac{\partial\varphi_j}{\partial x},y,
-\frac{\partial\varphi_j}{\partial y}\right): 
\frac{\partial\varphi_j}{\partial\theta_j}=0\right\}.$$
For $G\subset T^*Y$,
let $$\Lambda_j(G)=\{z\in T^*X: \exists\zeta\in G: 
(z,\zeta)\in\Lambda_j\}.$$
Then we have the following statement on the propagation 
of singularities, see \cite{Kamotski2007}.

\begin{thm}
Let $u\in {\mathcal{D}}^\prime(Y).$ Then
$WF(Iu)\subset \bigcup_j \Lambda_j(WF(u)).$
\label{th:singularities}
\end{thm}
Consequently, combining these observations with Theorem \ref{thm:SolRep} we obtain the following property.

\begin{cor} 
\label{thm:SolRepwf}
Let $n\geq 1$, $m\geq 2$, and let
\[
\left\{\begin{array}{l}
D_t u = A(x,D_x)u + B(t,x,D_x)u + f(t,x), \quad (t,x) \in [0,T] \times \R^n, \\
\left.u \right|_{t=0} = u_0(x), \quad x \in \R^n,
\end{array} \right.
\]
where $A(x,D_x)$ is an upper-triangular matrix of pseudo-diffe\-ren\-tial operators of order $1$ and $B(t,x,D_x)$ is a matrix of pseudo-differential operators of order $0$, continuous with respect to $t$. Recall that under condition (H1) and (H2), for any $N\in\N$, the components $u_j$, $j=1,\ldots,m$, of the solution $u$ are given by
	\begin{equation}
	\label{repr_formulawp}
	u_j(t,x) = \sum_{l=1}^{m} \left( \mc H_{j,l}^{l-j}(t) + R_{j,l}(t) \right) u_l^0 + \left( \mc K_{j,l}^{l-j}(t) + S_{j,l}(t) \right) f_l,
	\end{equation} 
	where $R_{j,l}$, $S_{j,l} \in \mc L(H^s, C([0,T],H^{s+N-l+j}))$ and the operators  
	$$\mc H_{j,l}^{l-j}, \mc K_{j,l}^{l-j}\in  \mc L(C([0,T], H^s), C([0,T],H^{s-l+j}))$$ are integrated Fourier Integral Operators of order $l-j$.
	
Consequently, up to any Sobolev order (depending on $N$), the wave front set of $u_j$ is given by 
\begin{equation}\label{EQ:wff}
WF(u_j(t,\cdot))\subset \left(\bigcup_{l=1}^m WF(\mc H_{j,l}^{l-j}(t) u_l^0)\right)
\bigcup
\left(\bigcup_{l=1}^m WF(\mc K_{j,l}^{l-j}(t) f_l)\right),
\end{equation} 
with each of the wave front sets for terms in the right hand side of \eqref{EQ:wff} given by the 
propagation along the broken Hamiltonian flow as in Theorem \ref{th:singularities}.
\end{cor}

We conclude the paper by presenting some applications of Theorem \ref{thm:SolRep} and Theorem \ref{main_theo_mi}.

\section{Application: Higher order hyperbolic equations}

In this section we want to study the well-posedness of the Cauchy problem
\begin{equation}
\label{CP}
\left\{
\begin{array}{cc}
D^m_t u=\lsum_{j=0}^{m-1} A_{m-j}(t,x,D_x)D_t^j u+f(t,x),&\quad (t,x)\in[0,T]\times\R^n,\\
D^{k-1}_t u(0,x)=g_{k}(x),&\quad k=1,...,m,
\end{array}
\right.
\end{equation}
where each $A_{m-j}(t,x,D_x)$ is a scalar differential operator of order $m-j$ with continuous and bounded coefficients depending on $t$ and $x$.
As usual, $D_t=\frac{1}{\irm}\partial_t$ and $D_x=\frac{1}{\irm}\partial_x$. Let $A_{(m-j)}$ denote the principal part of the operator $A_{m-j}$, We assume that the problem above is hyperbolic, i.e., the characteristic equation
\[
\tau^m=\sum_{j=0}^{m-1}A_{(m-j)}(t,x,\xi)\tau^j
\equiv\sum_{j=0}^{m-1}\sum_{|\gamma|=m-j}a_{m-j,\gamma}(t,x)\xi^\gamma\tau^j.
\]
has $m$ real valued roots $\lambda_1, \lambda_2,\cdots, \lambda_m$. In addition we work under the hypothesis that the roots $\lambda_i$, $i=1,\dots,n$ are symbols of order $1$, i.e.,
\[
\lambda_i\in C([0,T], S^1),
\]
for all $i=1,\dots,n$. For this reason we assume that
\begin{itemize}
\item[({\bf H0})] the coefficients of the equation above are continuous in $t$ and smooth in $x$, with bounded derivatives of any order $\alpha\in\N_0^n$ with respect to $x$. 
\end{itemize}
We will make use first of Theorem \ref{main_theo_mi} and then of Theorem \ref{thm:SolRep}.

\subsection{Well-posedness in Sobolev spaces}
We begin by reducing the $m$-order partial differential equation in \eqref{CP} into a first order system of pseudo-differential equations. Let  $\lara{D_x}$ be
the pseudo-differential operator with symbol $\lara{\xi}$. The transformation
\[
u_k=D_t^{k-1}\lara{D_x}^{m-k}u,
\]
with $k=1,...,m$, makes the Cauchy problem \eqref{CP} equivalent to the following system
\begin{equation}
\label{syst_Taylor}
D_t\left(
\begin{array}{c}
u_1 \\
\cdot \\
\cdot\\
u_m \\
\end{array}
\right)
= \left(
\begin{array}{ccccc}
0 & \lara{D_x} & 0 & \dots & 0\\
0 & 0 & \lara{D_x} & \dots & 0 \\
\dots & \dots & \dots & \dots & \lara{D_x} \\
b_1 & b_2 & \dots & \dots & b_m \\
\end{array}
\right)
\left(\begin{array}{c}
u_1 \\
\cdot \\
\cdot\\
u_m \\
\end{array}
\right)
+
\left(\begin{array}{c}
0 \\
0 \\
\cdot\\
f \\
\end{array}
\right),
\end{equation}
where
\[
b_j=A_{m-j+1}(t,x,D_x)\lara{D_x}^{j-m},
\]
with initial condition
\begin{equation}
\label{ic_Taylor}
u_k|_{t=0}=\lara{D_x}^{m-k}g_k,\qquad k=1,...,m.
\end{equation}
The matrix in \eqref{syst_Taylor} can be written as $A+B$ with
\[
A=\left(
\begin{array}{ccccc}
0 & \lara{D_x} & 0 & \dots & 0\\
0 & 0 & \lara{D_x} & \dots & 0 \\
\dots & \dots & \dots & \dots & \lara{D_x} \\
b_{(1)} & b_{(2)} & \dots & \dots & b_{(m)} \\
\end{array}
\right),
\]
where $b_{(j)}=A_{(m-j+1)}(t,x,D_x)\lara{D_x}^{j-m}$ and
\[
B=\left(
\begin{array}{ccccc}
0 & 0 & 0 & \dots & 0\\
0 & 0 & 0& \dots & 0 \\
\dots & \dots & \dots & \dots & 0 \\
b_1-b_{(1)} & b_2-b_{(2)} & \dots & \dots & b_m-b_{(m)} \\
\end{array}
\right).
\]
It is clear that the eigenvalues of the symbol matrix $A(t,x,\xi)$ are the roots $\lambda_j(t,\xi)$, $j=1,...,m$.

We want to apply Theorem \ref{main_theo_mi} to our Cauchy problem. This means to find under which hypotheses the equation in \eqref{CP} can be reduced into a first order system with upper-triangular principal part and lower order terms of suitable order as in (H1).

\begin{thm}
\label{theo_eq_1}
Let
	\[
	\left\{
	\begin{array}{ll}
	D^m_t u=\lsum_{j=0}^{m-1} A_{m-j}(t,x,D_x)D_t^j u+f(t,x),&\quad (t,x)\in[0,T]\times\R^n,\\
	D^{k-1}_t u(0,x)=g_{k}(x),&\quad k=1,...,m,
	\end{array}
	\right.
	\]
	where each $A_{m-j}(t,x,D_x)$ is a differential operator of order $m-j$ with continuous and bounded coefficients depending on $t$ and $x$ as in (H0). Let $A_{(m-j)}$ denote the principal part of the operator $A_{m-j}$. Assume that the roots of the corresponding characteristic polynomial are real valued symbols
	\[
	\lambda_i\in C([0,T], S^1), \quad 1,\dots,m,
	\] 
	and that  
	\[
	A_{m-j+1}(t,x,\xi)\in C([0,T], S^0)
	\]
	for all $j=1,\dots,m-1$. If $f\in C([0,T], H^{s+m-1})$ and $g_k\in H^{s+m-1}(\R^n)$ for all $k=1,\dots,m$ then the Cauchy problem \eqref{CP} has a unique solution $u\in C^{m-1}([0,T], H^{s+m-1})$.
\end{thm}
\begin{proof}
	We consider the associated reduced system with principal part given by the matrix 
	\[
	A=\left(
	\begin{array}{ccccc}
	0 & \lara{D_x} & 0 & \dots & 0\\
	0 & 0 & \lara{D_x} & \dots & 0 \\
	\dots & \dots & \dots & \dots & \lara{D_x} \\
	b_{(1)} & b_{(2)} & \dots & \dots & b_{(m)} \\
	\end{array}
	\right),
	\]
	where $b_{(j)}=A_{(m-j+1)}(t,x,D_x)\lara{D_x}^{j-m}$. The operators $b_{(j)}$ are of order 1 so since we assume that $A_{m-j+1}(t,x,\xi)\in C([0,T], S^0)$ for $j=1,\dots,m-1$ it follows that $b_{(j)}\equiv 0$ for $j=1,\dots,m-1$. This means that the Sylvester matrix $A$ is actually upper-triangular and that the matrix $B$ of the lower order terms is of the following type:
	\[
	B=\left(
	\begin{array}{ccccc}
	0 & 0 & 0 & \dots & 0\\
	0 & 0 & 0& \dots & 0 \\
	\dots & \dots & \dots & \dots & 0 \\
	b_1 & b_2 & \dots & \dots & b_m-b_{(m)} \\
	\end{array}
	\right),
	\]
	with $b_j=A_{m-j+1}(t,x,D_x)\lara{D_x}^{j-m}$, for $j=1,\dots,m-1$ and $b_m-b_{(m)}=A_1(t,x,D_x)-A_{(1)}(t,x,D_x)$. Since  $A_{m-j+1}(t,x,\xi)\in C([0,T], S^0)$ for $j=1,\dots,m-1$ we have that $b_j$ is a pseudo-differential operator of order $j-m$ for $j=1,\dots,m$.
	We are therefore under the assumptions of Theorem \ref{main_theo_mi} for the matrices $A$ and $B$. Since $g_k\in H^{s+m-1}(\R^n)$ the initial data    
	\[
	\lara{D_x}^{m-k}g_k
	\]
	of the reduced Cauchy problem belong to the space $H^{s+k-1}(\R^n)$ for all $k=1,\dots,m$. Thus, by Theorem  \ref{main_theo_mi} there exists a unique solution $u(t,x)$ to the Cauchy problem under consideration such that
	\[
	D_t^{k-1}\lara{D_x}^{m-k}u\in C([0,T], H^{s+k-1}),
	\]
	for $k=1,\dots,m$. By Sobolev mapping properties of pseudo-differential operators it follows that $u\in C^{m-1}([0,T], H^{s+m-1})$.
\end{proof}

\subsubsection{Second order examples}
\label{subsec_ex_1}

\begin{itemize}
\item[(i)]
The equation
\[
D_t^2u=a_2(t,x)u+(a_1(t,x)D_x+b(t,x))D_tu+f(t,x),
\]
where $x\in\R$ and $a_1$ is real valued falls into the class of equations considered in the previous theorem. Indeed, the characteristic polynomial
\[
\tau^2-a_1(t,x)\tau\xi
\]
has two real roots and $A_2=a_2(t,x)$ is an operator of order $0$.

\item[(ii)] Let us now consider the second order Cauchy problem
\[
	\left\{
	\begin{array}{ll}
	&D^2_t u=a^2(t)D_x^2 u +b_1(t)D_xu+b_2(t)D_tu+b_3(t)u+f(t,x), \\
	&u(0,x)=g_{0}(x),\\
	&D_tu(0,x)=g_1(x),\\
	\end{array}
	\right.
	\]
	where $(t,x)\in[0,T]\times\R$, the equation coefficients are continuous and $a\in C^1$ with $a(t)\ge 0$. Making use of the standard reduction into first order system of pseudo-differential equations we have that the Cauchy problem above is equivalent to 
	\begin{equation}
	\label{syst_example_t}
	\begin{split}
	D_tU&=\left(
	\begin{array}{cc}
	0 & \lara{D_x}\\ 
	a^2(t)D_x^2\lara{D_x}^{-1} & 0
	\end{array}
	\right)U
	\\
	&+\left(
	\begin{array}{cc}
	0 & 0\\ 
	b_1(t)D_x\lara{D_x}^{-1}+b_3(t)\lara{D_x}^{-1} & b_2(t)
	\end{array}
	\right)U  
	+\left(
	\begin{array}{c}
	0\\ 
	f(t,x)
	\end{array}
	\right),
	\end{split}
	\end{equation}
	where $U=(u_1,u_2)^T=(\lara{D_x}u, D_tu)^T$ and $U(0,x)=U_0=(\lara{D_x}g_0,g_1)^T$. In the sequel we will denote the right-hand side of the system above with 
	\[
	A(t,D_x)U+B(t,D_x)+F,
	\]
	where $A$ and $B$ are defined by operators of order $1$ and $0$, respectively. The principal part matrix 
	\[
	A(x,\xi)=\left(
	\begin{array}{cc}
	0 & \lara{\xi}\\ 
	a^2(t)\xi^2\lara{\xi}^{-1} & 0
	\end{array}
	\right)
	\]
	is not upper triangular and it is not diagonalisable because of the zeros of the coefficient $a$. However it can be reduced into upper-triangular form. We refer here the reader to \cite{Gramchev2013} and to the appendix in \cite{Garetto2018}. The matrix $A$ has $\lambda_1(t,\xi)=-a(t)\xi$ and $\lambda_2(t,\xi)=a(t)\xi$ as eigenvalues. It follows that 
	\[
	h^{(1)}=\left(
	\begin{array}{c}
	1\\
	a(t)\xi\lara{\xi}^{-1}\\
	\end{array}
	\right)
	\]
	and
\[
	h^{(2)}=\left(
	\begin{array}{c}
	1\\
	-a(t)\xi\lara{\xi}^{-1}\\
	\end{array}
	\right)
	\]	
	are eigenvectors corresponding to $\lambda_1$ and $\lambda_2$, respectively. Choose now $h^{(1)}$ (the argument is analogous with $h^{(2)}$). The matrix
	\[
	T_1=(h^{(1)}, e_2)=\left(
	\begin{array}{cc}
	1 & 0\\ 
	a(t)\xi\lara{\xi}^{-1} & 1
	\end{array}
	\right)
	\]
	is invertible. Its inverse is 
	\[
	T_1^{-1}=\left(
	\begin{array}{cc}
	1 & 0\\ 
	-a(t)\xi\lara{\xi}^{-1} & 1
	\end{array}
	\right)
	\]
	and
	\[
	T_1^{-1}AT_1=\left(
	\begin{array}{cc}
	-a(t)\xi & \lara{\xi}\\ 
             0 & a(t)\xi
	\end{array}
	\right)
	= \left(
	\begin{array}{cc}
	\lambda_1(x,\xi) & \lara{\xi}\\ 
             0 & \lambda_2(x,\xi)
	\end{array}
	\right).
	\]
	Note that the operator $T_1^{-1}(t,D_x)A(t,D_x)T_1(t,D_x)$ can be therefore written as 
	\[
	\left(\begin{array}{cc}
	\lambda_1(t,D_x) & \lara{D_x}\\ 
             0 & \lambda_2(t,D_x)
	\end{array}
	\right). 	
	\]
 We can now use this transformation to reduce the system $D_tU=AU$ into upper-triangular form. More precisely, for $U=T_1V$, we have that the system
 \[
 D_tU=A(t,D_x)U+B(t,D_x)U+F
 \] 
is equivalent to 
\[
\begin{split}
D_tV&=(T_1^{-1}AT_1)(t,D_x)V+T_1^{-1}(BT_1+D_tT_1)V+T_{1}^{-1}F,\\
&=\left(\begin{array}{cc}
	\lambda_1(t,D_x) & \lara{D_x}\\ 
             0 & \lambda_2(t,D_x)
	\end{array}
	\right)V+ \left(\begin{array}{cc}
	0 & 0\\ 
         b(t,\xi) & b_2(t)
	\end{array}
	\right)V+F,
\end{split}
\]
where
\[
b(t,\xi)=(b_1(t)+b_2(t)a(t)+D_t a(t))\xi\lara{\xi}^{-1}+b_3(t)\lara{\xi}^{-1}
\]

 	%\[
	%\begin{split}D_tV=&\left(
	%\begin{array}{cc}
	%\lambda_1(x,D_x) & \lara{D_x}\\ 
      %       0 & \lambda_2(x,D_x)
	%\end{array}
%	\right)V+\left(\begin{array}{cc}
%	0& 0\\ 
    %         b_{2,1}(x,D_x) & 0
%	\end{array}
%	\right)V+F,\\
	%V(0,x)=&T^{-1}U_0.
%	\end{split}
%	\]
 
Under the assumptions that $a\in C^1([0,T])$  we easily see that the eigenvalues $\lambda_1$ and $\lambda_2$ belong to $C([0,T], S^1)$. In addition, condition (H1) is fulfilled if the symbol above is of order $-1$. This is the case when 
\[
b_1(t)+b_2(t)a(t)+D_t a(t)=0,
\]
for all $t\in[0,T]$ (for instance when $b_2\equiv 0$ and $b_1=-D_ta$).
Since $T^{-1}$ is a matrix of pseudo-differential operators of order $0$ we have that $V(0,x)=V_0=T_{1}^{-1}(0,D_x)U_0$ has the same regularity properties of $U_0$.  We can therefore apply Theorem \ref{thm:SolRep} to this system and obtain the following result.

\begin{thm}
\label{theo_wave_eq}
Let 
\[
	\left\{
	\begin{array}{ll}
	&D^2_t u=a^2(t)D_x^2 u +b_1(t)D_xu+b_2(t)D_tu+b_3(t)u+f(t,x), \\
	&u(0,x)=g_{0}(x),\\
	&D_tu(0,x)=g_1(x),\\
	\end{array}
	\right.
\]
	where $a(t)\ge 0$ is of class $C^1$, the lower order coefficients $b_i(t)$, $i=1,2,3$, are continuous and
\[
b_1(t)+b_2(t)a(t)+D_t a(t)=0,
\]
for all $t\in[0,T]$. Let $s\in\R$. If $f\in C([0,T], H^{s+1})$ and $g_0,g_1\in H^{s+1}(\R^n)$ then there exists a unique solution $u\in C^1([0,T], H^{s+1})$ of the Cauchy problem above.
\end{thm}
In particular it follows that if $f\in C^\infty([0,T]\times\R)$ and compactly supported with respect to $x$ and $g_0, g_1\in C^\infty_c(\R)$ then this Cauchy problem is $C^\infty$ well-posed. 
%This result is consistent with the well-posedness result in \cite{SpaTaglia07} where symmetrisation techniques have been used rather than Fourier integral operator techniques. Note also that since the coefficient $a$ depends only on $x$ the Oleinik condition (see \cite{Oleinik70}) is trivially fulfilled. This guarantees the $C^\infty$ well-posedness of the Cauchy problem above. 

 \end{itemize}

\subsection{Representation formula for the solution}
We now assume that the principal part of the equation 
\[
D^m_t u=\lsum_{j=0}^{m-1} A_{m-j}(t,x,D_x)D_t^j u+f(t,x)
\]
depends only on $x$, i.e.,
\[
D^m_t u=\lsum_{j=0}^{m-1} A_{(m-j)}(x,D_x)D_t^j u+\sum_{j=0}^{m-1} (A_{m-j}-A_{(m-j)})(t,x,D_x)D_t^j u+f(t,x).
\]
If in addition to the hypotheses of Theorem \ref{theo_eq_1} we assume that the roots of this equation fulfil condition (H2) we obtain the following representation formula by straightforward application of Theorem \ref{thm:SolRep}.
\begin{thm}
\label{theo_eq_2}
Let
	\[
	\left\{
	\begin{array}{ll}
	D^m_t u=\lsum_{j=0}^{m-1} A_{m-j}(t,x,D_x)D_t^j u+f(t,x),&\quad (t,x)\in[0,T]\times\R^n,\\
	D^{k-1}_t u(0,x)=g_{k}(x),&\quad k=1,...,m,
	\end{array}
	\right.
	\]
	where each $A_{m-j}(t,x,D_x)$ is a differential operator of order $m-j$ with
	with continuous and bounded coefficients depending on $t$ and $x$ as in (H0). Let the principal part $A_{(m-j)}$ of the operator $A_{m-j}$ be independent of $t$ and let the roots $\lambda_i$, $i=1,\dots, m$, of the corresponding characteristic polynomial be real valued symbols of order $1$ fulfilling condition (H2). Assume that 
 \[
	A_{m-j+1}(t,x,\xi)\in C([0,T], S^0)
\]
	for all $j=1,\dots,m-1$. If $f\in C([0,T], H^{s+m-1})$ and $g_k\in H^{s+m-1}(\R^n)$ for all $k=1,\dots,m$, then for any $N\in\N$ and $j=1,\dots,m$ the solution $u\in C^{m-1}([0,T], H^{s+m-1})$ to the Cauchy problem can be written as    
\begin{equation}
	\label{repr_formula_eq}
	D_t^{j-1}u(t,x) = \sum_{l=1}^{m} \left( \mc H_{j,l}^{l-m}(t) + R_{j,l}(t) \right) \lara{D_x}^{m-l}g_l + \left( \mc K_{j,m}^{0}(t) + S_{j,m}(t) \right) f, 
		\end{equation} 
	where 
	\begin{itemize}
	\item[(i)]  $\mc H_{j,l}^{l-m}$ is an integrated Fourier Integral Operator of order $l-m$,
	\item[(ii)] $\mc K_{j,m}^{0}$ is an  integrated Fourier Integral Operator of order $0$,
	\item[(iii)] $R_{j,l}\in \mc L(H^s, C([0,T],H^{s+N-l+m}))$,
	\item[(iv)] $S_{j,m} \in \mc L(H^s, C([0,T],H^{s+N}))$.
 
	\end{itemize}

\end{thm} 
\begin{proof}
Since we work under the assumptions (H1) and (H2) we can apply Theorem \ref{thm:SolRep} to the system in \eqref{syst_Taylor} where the matrix $A$ of the principal part is only depending on $x$ and is upper-triangular. Note that the right-hand side is of the
type  
\[
\left(\begin{array}{c}
0 \\
0 \\
\vdots\\
f \\
\end{array}
\right)
\]
and the initial condition is given by 
\[
u_k|_{t=0}=\lara{D_x}^{m-k}g_k,\qquad k=1,...,m.
\]
A straightforward application of Theorem \ref{thm:SolRep} allows us to write
\[
 D_j^{j-1}\lara{D_x}^{m-j}u=\sum_{l=1}^{m} \left( \mc H_{j,l}^{l-j}(t) + R_{j,l}(t) \right) \lara{D_x}^{m-l}g_l + \left( \mc K_{j,m}^{m-j}(t) + S_{j,m}(t) \right) f,
\]
where $R_{j,l} \in \mc L(H^s, C([0,T],H^{s+N-l+j}))$,  $S_{j,m} \in \mc L(H^s, C([0,T],H^{s+N-m+j}))$ and the operators  $\mc H_{j,l}^{l-j}$ and $\mc K_{j,m}^{m-j}$ are of order $l-j$ and $m-j$, respectively. 
	It follows that
	\begin{multline*}
 D_j^{j-1}u=\sum_{l=1}^{m}\lara{D_x}^{-m+j}\left( \mc H_{j,l}^{l-j}(t) + R_{j,l}(t) \right) \lara{D_x}^{m-l}g_l \\
 + \lara{D_x}^{-m+j}\left( \mc K_{j,m}^{m-j}(t) + S_{j,m}(t) \right) f,
\end{multline*}
for all $j=1,\dots,m$. By composition of pseudo- and Fourier integral operators we easily see that 
\begin{itemize}
\item[-] $\lara{D_x}^{-m+j}\mc H_{j,l}^{l-j}(t)$ is of order $l-m$,
\item[-] $\lara{D_x}^{-m+j}R_{j,l} \in \mc L(H^s, C([0,T],H^{s+N-l+m}))$,
\item[-] $\lara{D_x}^{-m+j}\mc K_{j,m}^{m-j}$ is of order $0$,
\item[-]  $\lara{D_x}^{-m+j}S_{j,m}\in \mc L(H^s, C([0,T],H^{s+N}))$.
\end{itemize}
This shows that $u$ can be written as in \eqref{repr_formula_eq} and completes the proof.
	
 \end{proof}
 
 \subsubsection{Example}
 Let us consider an $m$ order homogeneous differential operator
 \[
 A(x,D_t,D_x),\qquad t\in[0,T],\, x\in\R,
 \]
 such that its symbol (which coincides with its principal symbol) is
 \[
 A(x,\tau,\xi)=\Pi_{i=1}^m (\tau-a_i(x)\xi).
 \]
Assume that all the coefficients $a_i$ are real valued, smooth and bounded. Assume that the derivatives of the coefficients $a_i$, $i=1,\dots,n$, are bounded as well. It follows that the roots of the characteristic polynomials above are symbols of order 1 and the operator $A$ is in general a hyperbolic operator with multiplicities. We work under the assumption that when the equation $A(x,D_t,D_x)u=f$ is transformed into a first order system of pseudo-differential equations the matrix of the system is upper-triangular (for a second order example see Subsection \ref{subsec_ex_1}(i)). Since there are no lower order terms then condition (H1) is trivially fulfilled.
%Note that when $m=2$ the operator $A$ is of the form
%\[
%D_t^2-(a_1(x)+a_2(x))D_tD_x+a_1(x)a_2(x)D_x^2.
%\]
%Since the coefficients $a_i$, $i=1,\dots,m$ are real valued this is in general a hyperbolic operator with multiplicities. Arguing as in paragraph (ii) in the previous subsection we have that the well-posedness of the Cauchy problem
We can now investigate the well-posedness of the Cauchy problem
\begin{equation}
\label{CP-example_m}
\begin{split}
A(x,D_t,D_x)u&=f(t,x),\\
D_t^{k-1}u(0,x)&=g_k,
\end{split}
\end{equation}
$k=1,\dots,m-1$. We have that if $f\in C^\infty([0,T],\R)$ and $g_k\in C^\infty_c(\R)$ for all $k=1,\dots,m-1$ the Cauchy problem \eqref{CP-example_m} has a unique solution $u\in C^\infty([0,T]\times\R)$. In addition if the roots 
\[
\lambda_i(x,\xi)=a_i(x)\xi,\quad i=1,\dots,m
\]
fulfil the condition (H2) then the representation formula \eqref{repr_formula_eq} holds. This for instance happens when the coefficients $a_i$ have distinct first derivatives in the points of multiplicities, i.e. $a_i(x)=a_j(x)$ implies $a_i'(x)\neq a_j'(x)$.

\bibliographystyle{alphaabbr}
\bibliography{LitMicrolocal}

\end{document}